\documentclass[11pt]{article}%
\usepackage{amsmath,amsthm,amssymb,amscd}
\usepackage{latexsym}
\usepackage{epsfig}
\usepackage{graphics}
\usepackage{color}
\usepackage{graphicx}
\usepackage{amsmath}
\usepackage{amsfonts}
\usepackage{amssymb}%
\providecommand{\U}[1]{\protect\rule{.1in}{.1in}}
\renewcommand{\baselinestretch}{1.5}
\makeatletter \@addtoreset{equation}{section}

\makeatother

\def\singlespace{\def\baselinestretch{1}\@normalsize}

\def\T{{ \mathrm{\scriptscriptstyle T} }}

\newcommand{\Date}[1]{\def\@Date{#1}}
\def\today{\number\day~\ifcase\month\or
 January\or February\or March\or April\or May\or June\or
 July\or August\or September\or October\or November\or December\fi~\number\year}

\begin{document}

\title{\bf A Double AR Model Without Intercept: an Alternative to Modeling Nonstationarity and Heteroscedasticity }
\author{Dong Li$^1$, Shaojun Guo$^{2,3}$ and Ke Zhu$^{2}$\\
Tsinghua University$^1$, Chinese Academy of Sciences$^2$, and London School of Economics$^3$\\ {dongli@math.tsinghua.edu.cn}\ \ {guoshaoj@amss.ac.cn} \ \ {kzhu@amss.ac.cn}}

\maketitle
\date{}

\begin{abstract}
This paper presents a \underline{d}ouble \underline{AR} model \underline{w}ithout \underline{in}tercept (DARWIN
model) and provides us a new way to study the non-stationary heteroskedastic time
series. It is shown that the DARWIN model is always non-stationary and heteroskedastic, and its sample properties depends on the Lyapunov exponent. An
easy-to-implement estimator is proposed for the Lyapunov exponent, and it is unbiased, strongly consistent and asymptotically normal.
Based on this estimator, a powerful test is constructed for testing the
stability of the model. Moreover, this paper proposes the
quasi-maximum likelihood estimator (QMLE) for the DARWIN model, which has an explicit form. The strong consistency and asymptotical normality of the
QMLE are established regardless of the sign of the Lyapunov
exponent. Simulation studies are conducted to assess the performance
of the estimation and testing and an empirical example is given for
illustrating the usefulness of the DARWIN model.
\end{abstract}

%\begin{keywords}

\newpage

\section{Introduction}

Time-varying volatility has been crucial in modeling economic and financial time series.
Since the seminal work of Engle (1982), the autoregressive
conditional heteroscedasticity (ARCH) model and its numerous
variants have been widely used; see, e.g., Bollerslev {\it et al.} (1992) and Francq and Zako\"{i}an (2010). Among them, the first order double autoregressive (DAR)
model has attracted much attention, which takes the form
\begin{equation}
\label{dar}
y_t=\phi_0 y_{t-1}+ \eta_t\sqrt{\omega_0+\alpha_0 y_{t-1}^2},
\end{equation}
where $\phi_0\in (-\infty,\infty)$, $\omega_0 \ge 0$, $\alpha_0\geq0$,
$\{\eta_t\}$ is a sequence of independent and identically
distributed (i.i.d.) random variables with zero mean and unit
variance, and $\eta_t$ is independent of $\{y_{j}; j<t\}$. Model
(\ref{dar}) was initially introduced by Weiss (1984), and the term
`\textit{double autoregressive}' was coined by Ling (2004) since
both the conditional mean and variance functions are regressions
only on the observed data and its square, respectively. Clearly,
model (\ref{dar}) belongs to the class of \textsc{ARMA-ARCH} models
in Weiss (1984) and of  nonlinear AR models in Cline and Pu (2004), but it
is different from Engle's \textsc{ARCH} model if $\phi_0\neq 0$. Its higher-order extension and generalization can be found in Weiss (1984), Lu (1998), Ling (2007), Zhu and Ling (2013), Guo, Ling and Zhu (2014), Li, Ling and Zako\"{i}an (2015), Li, Ling and Zhang (2015) and many others.

The stationarity conditions and the associate inferential theory of model (\ref{dar}) have been well studied under the compact parameter space $\Lambda = \{\lambda_{0} = (\phi_{0},\omega_{0},\alpha_{0}): |\phi_{0}| \le \bar{\phi}, \underline{\omega} \le \omega_0 \le \bar{\omega}, \underline{\alpha} \le \alpha_0 \le \bar{\alpha}\}$, where $\bar{\phi}$, $\underline{\omega}$, $\bar{\omega}$, $\underline{\alpha}$ and $\bar{\alpha}$ are some finite positive constants.
Early contributions in this context include Gu$\acute{e}$gan and Deibolt (1994) and Borkovec and Kl$\ddot{u}$ppelberg (2001), who derived the sufficient and necessary condition of weak stationarity of model (\ref{dar}). Recently, Chen, Li and Ling (2014) proved that when $\eta_t$ is symmetric, model (\ref{dar}) is strictly stationary if and only if the Lyapunov exponent $\gamma_0=E\log|\phi_0 +\eta_t\sqrt{\alpha_0}|<0$.
%Figure \ref{boundary} plots the stationary region of model (\ref{dar}) for three different distributions of $\eta_{t}$. %From this figure, it is observed that, different from the traditional AR(1) model,  model (\ref{dar}) is stationary even if $|\phi_{0}|$ is 1 or larger than 1.

By assuming that the true value $\lambda_0$ is an interior point of $\Lambda$, the
inference of model (\ref{dar}) based on quasi-maximum likelihood estimation exhibits quite different phenomenon
in terms of $\gamma_{0}$. For instance, when $\gamma_{0}<0$, Ling (2004) showed that the quasi-maximum likelihood estimator (QMLE) of $\lambda_{0}$ is consistent and asymptotically
normal; and when $\gamma_{0}\geq0$, Ling and Li(2008) and Chen, Li and Ling (2014)
demonstrated that the (unconstrained) QMLE of $(\phi_0, \alpha_0)$
is consistent and asymptotically normal, but the intercept term $\omega_{0}$ cannot be consistently estimated.
This phenomenon can also be found for the least absolute deviation estimator and robust quasi-maximum likelihood estimation in Chan and Peng (2005), Zhu and Ling (2013), and Li, Guo, and Li (2015).

%However, all the above results require the innovation error to have a finite fourth moment. To relax this, \cite{cp}, \cite{zhu}, and \cite{lgl} proposed the least absolute deviation estimation and robust quasi-maximum likelihood estimation, respectively. Moreover, \cite{lgl} provided a consistent estimator of $\gamma_0$ without stationarity assumption and proposed a robust procedure for strict stationary testing.

In all aforementioned work, the positivity of $\omega_0$ is essential for the strict stationarity and (robust) quasi-maximum likelihood estimation, although the intercept $\omega_0$ is not involved directly in $\gamma_0$. The case that $\omega_0 = 0$ would be meaningful but hardly touched so far.
As one motivation of this paper, it is of interest to fill in this gap
from a theoretical viewpoint. Another parallel motivation of this study is from application. The
importance of model (\ref{dar}) and its higher-order extension has
been well demonstrated by empirical studies; see, e.g., Ling (2004), Zhu and Ling (2013), and many others. In these applications, however, an interesting
finding is that the intercept term $\omega_0$ is often very close to
zero; see also the illustrating example in Section 6 below. As a result, it is intuitively appealing to study a first-order \underline{d}ouble \underline{AR} model \underline{w}ithout \underline{in}tercept (abbreviated to
DARWIN(1)) as follows:
\begin{equation}\label{model}
y_t=\phi_0 y_{t-1}+\eta_t\sqrt{\alpha_0 y_{t-1}^2},
\end{equation}
where all notations inherit from model (\ref{dar}) with
$\alpha_{0}>0$. In the special case of $\phi_0=0$, model (\ref{model})
becomes the ARCH model without intercept in Hafner and Preminger (2015). 
In this paper, we are concerned with the probability structure of model (\ref{model}) and the estimation and inference of $\gamma_0$ and $(\phi_0,\alpha_0)$.
It is surprising that its probabilistic structures and asymptotics
of the QMLE are totally different from those of model (\ref{dar}).
This is out of our expectation. Moreover, model (\ref{model}) is always  non-stationary and heteroskedastic regardless of the sign of $\gamma_{0}$, and hence it provides us a new way to model the non-stationary heteroskedastic time series. Specifically, it can be seen that when $\gamma_0=0$, the conditional volatility of $y_{t}$ is a
nondegenerate random variable oscillating randomly between zero and infinity over time, while when
$\gamma_{0}\not=0$, this is not the case. In this sense, model (2) is stable if $\gamma_{0}=0$, and unstable otherwise; see also Hafner and Preminger (2015) for the same argument in the ARCH model without intercept.

This paper is organized as follows. Section 2
considers sample path properties of $y_{t}$ in model (\ref{model}). Section 3
proposes a new estimator of $\gamma_0$ and discusses its asymptotic
theory and also a test for stability. Section 4 focuses on the
QMLE of $(\phi_{0},\alpha_{0})$ in model (\ref{model}) and derives its asymptotic properties. Sections 5 and 6 investigate the numerical properties of the proposed procedures using both simulated and real data. The conclusions are
offered in Section 7. All technical proofs are relegated to the
Appendix.

\section{Sample path properties}
%Let $\eta$ be a generic random variable with the same distribution
%as $\eta_t$.
In this section, we study sample path properties of $y_{t}$ in model
(\ref{model}), when $\eta_{t}$ is symmetric but not necessarily has
mean zero and variance one. To accomplish it, we consider an
auxiliary process
\begin{equation}\label{auxiliary}
x_t=\Big|\phi_0x_{t-1}+\eta_t\sqrt{\alpha_0x^2_{t-1}}\Big|,
\end{equation}
where the initial value $x_0=|y_0|$. It is
straightforward to see that
\begin{equation}\label{x_process}
x_t=|\phi_0+\eta_t\sqrt{\alpha_0}|x_{t-1}\quad\mbox{or}\quad \log
x_t=\log|\phi_0+\eta_t\sqrt{\alpha_0}|+\log x_{t-1}.
\end{equation}
Note that when $\eta_{t}$ is symmetric, it is readily seen that
\begin{equation}\label{equiv_distribution}
\{x_t\}\stackrel{d}{=}\{|y_t|\};
\end{equation}
see, e.g., Borkovec (2000) and Borkovec and Kl$\ddot{u}$ppelberg (2001). Thus, $|y_t|$
either converges to zero or diverges to infinity almost surely at an exponential
rate as $t\rightarrow\infty$, according to $\gamma_0<0$ or
$\gamma_0>0$, respectively. This result sharply differs from the one
in model (\ref{dar}) when $\gamma_0<0$.

Next, we precisely characterize asymptotic
distribution of $|y_t|$ after a suitable renormalization.  Let $[a]$ be
the integral part of any real number $a$. The expression (\ref{x_process}) implies that for any $s\in[0, 1]$,
\begin{eqnarray*}
\frac{1}{\sqrt{n}}\log \frac{x_{[ns]}}{\exp\{[ns]\gamma_0\}}
=\frac{1}{\sqrt{n}}\sum_{i=1}^{[ns]}(\log|\phi_0+\eta_i\sqrt{\alpha_0}|-\gamma_0)
+\frac{1}{\sqrt{n}}\log x_0.
\end{eqnarray*}
By (\ref{equiv_distribution}) and Donsker's Theorem in Billingsley (1999, p.90), we have the following theorem:

{\bf Theorem 2.1.}
{\it
Suppose that $\{\eta_t\}$ is a sequence of i.i.d. symmetric random
variables satisfying
$\sigma_{\gamma_{0}}^2=\mathrm{var}(\log|\phi_0+\eta_{t}\sqrt{\alpha_0}|)\in(0,
\infty)$. If $y_0$ is symmetric with $P(y_0=0)=0$, and independent
of $\{\eta_t: t\geq 1\}$, then
\begin{eqnarray*}
\frac{|y_{{[ns]}}|^{\frac{1}{\sqrt{n}}}}{\exp(s\gamma_0\sqrt{n})}\Longrightarrow\exp\{\sigma_{\gamma_{0}}
\mathbb{B}(s)\} \quad \mbox{as $n\rightarrow\infty$ \quad in
$\mathbb{D}[0, 1]$},
\end{eqnarray*}
where `$\Longrightarrow$' denotes weak convergence, $\mathbb{B}(s)$ is a standard Brownian motion on $[0, 1]$,
and $\mathbb{D}[0, 1]$ is the space of functions defined on $[0, 1]$,
which are right continuous and have left limits, endowed with the Skorokhod topology.
%\end{theorem}
}

\begin{figure}[htp]
\begin{center}
% The arguments in the next line are {height}{optional width}{used only by OUP for typesetting}[filename, in directory art]
%\figurebox{15pc}{35pc}{}[sample_path.eps]
\includegraphics[height= 60mm,width= 140mm]{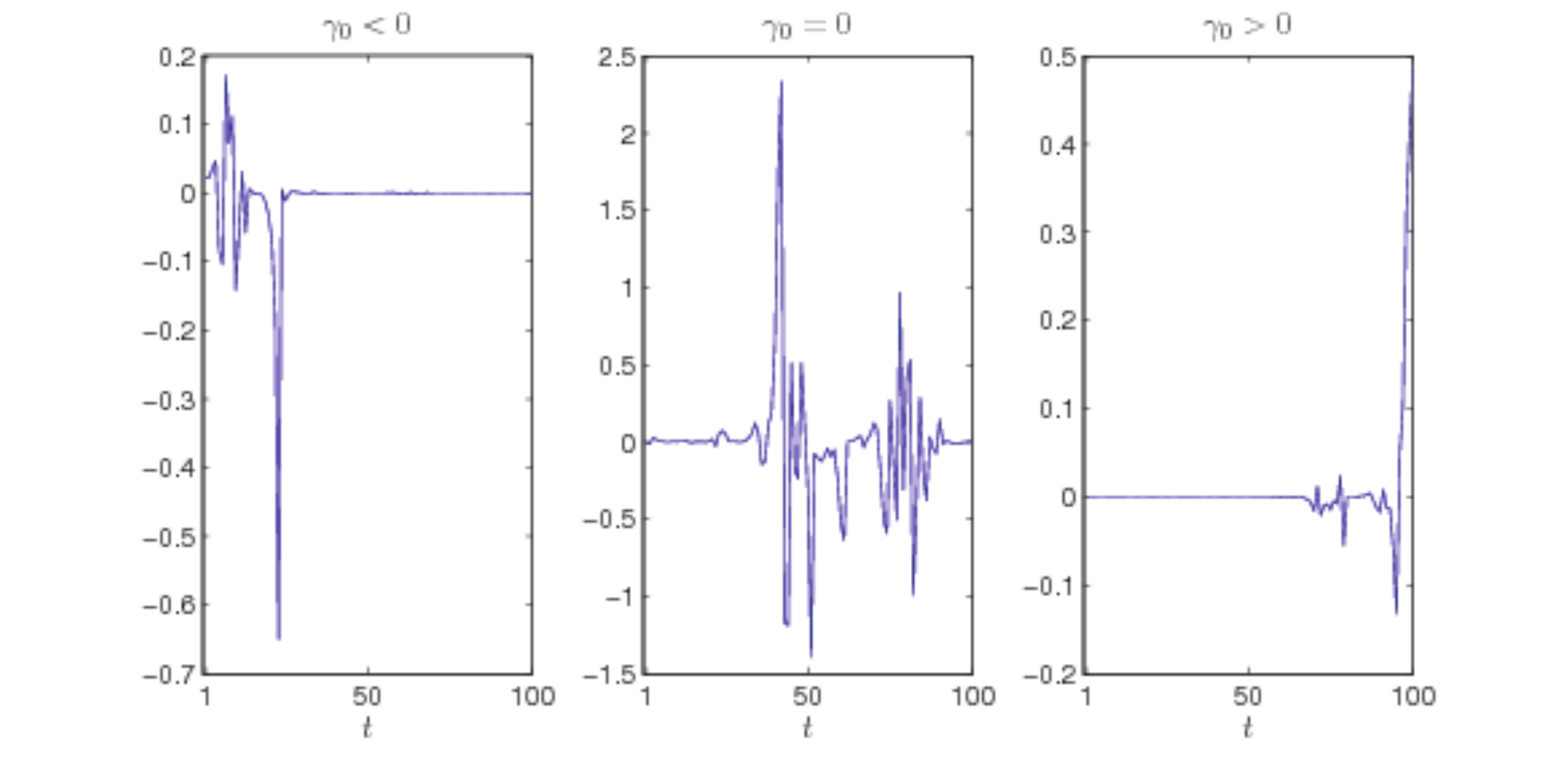}
% note that files may not be rotated
\caption{One sample path of $\{y_{t}\}_{t=1}^{100}$ in model (\ref{model}) corresponding to $\gamma_{0}<0$, $\gamma_{0}=0$, and $\gamma_{0}>0$, respectively.}
\label{sample_path}
\end{center}
\end{figure}

Theorem 2.1 has many implications. First, it implies that $y_{t}$ is non-stationary and heteroskedastic
regardless of the sign of $\gamma_{0}$. Second, its sample path property
depends on the sign of $\gamma_0$. Precisely, it indicates that, if $\gamma_0=0$,  the conditional volatility of $y_{t}$ is a nondegenerate random variable oscillating randomly between zero and infinity over time; otherwise, $|y_t|$  either converges to zero or diverges to infinity almost surely (a.s.) as $t \to \infty$. To illustrate this,  we present one sample path of $\{y_{t}\}$ of model (\ref{model}) with $\eta_{t}\sim N(0,1)$, $\phi_{0}=0.5$, and $\alpha_{0}=3.1$
(i.e., $\gamma_{0}<0$), 3.3058 (i.e., $\gamma_{0}=0$), or 3.5 (i.e.,
$\gamma_{0}>0$). The plots are depicted in Fig.\ref{sample_path}. They show clear evidence of different sample path properties with various $\gamma_0$. To distinguish them, model (\ref{model}) is called \textit{stable} if $\gamma_0=0$, and \textit{unstable} otherwise;
see Hafner and Preminger (2015) for the same argument in the ARCH model without
intercept. Since $\gamma_0$ plays a key role in determining the stability of model (\ref{model}), it is desirable to consider its estimation and inference in the next section.

\section{Estimation of Lyapunov exponent}

In this section, we propose a simple estimator for $\gamma_0$. This estimator
requires only the data $y_{n}$ and $y_0$ but shares good properties.
As in Section 2, we assume that $\eta_{t}$
is symmetric and not necessarily has mean zero and variance one. The basic idea is very intuitive. By (\ref{auxiliary}), we have
$\log(x_t/x_{t-1})=\log|\phi_0+\eta_t\sqrt{\alpha_0}|$, and then
\begin{eqnarray*}
\frac{1}{n}(\log x_n-\log
x_0)=\frac{1}{n}\sum_{t=1}^n\log(x_t/x_{t-1})=\frac{1}{n}\sum_{t=1}^n\log|\phi_0+\eta_t\sqrt{\alpha_0}|\to\gamma_{0}\,\,\,
\mbox{ a.s.}
\end{eqnarray*}
provided that $\gamma_{0}<\infty$. Consequently, an easy-to-implement
estimator of $\gamma_0$ is defined by
\begin{eqnarray*}
\widehat{\gamma}_n=\frac{1}{n}(\log |y_n|-\log |y_0|).
\end{eqnarray*}
The following theorem states that $\widehat{\gamma}_n$ is unbiased and
asymptotically normal.

{\bf Theorem 3.1.}
{\it
If the assumptions in Theorem 2.1 hold,  then,
$\mathrm{(i)}$ $\widehat{\gamma}_n$ is unbiased; $\mathrm{(ii)}$
$\widehat{\gamma}_n\rightarrow\gamma_0$ a.s.; and $\mathrm{(iii)}$
$\sqrt{n}(\widehat{\gamma}_n-\gamma_0)\Longrightarrow N(0,
\sigma^2_{\gamma_{0}})$ as $n\rightarrow\infty$, where
$\sigma^{2}_{\gamma_{0}}$ is defined as in Theorem 2.1.
}
%\end{theorem}

As an application, consider a testing problem whether model (\ref{model}) is stable or not, i.e.,
\begin{eqnarray*}
H_0: \gamma_0=0\quad \mbox{v.s.}\quad H_1: \gamma_0\neq 0.
\end{eqnarray*}
Based on Theorem 3.1, the proposed test statistic is
given by
\begin{equation}\label{test}
T_n=\sqrt{n}\frac{\widehat{\gamma}_n}{\widehat{\sigma}_\gamma},
\end{equation}
where
$\widehat{\sigma}^2_\gamma=\tfrac{1}{n}\sum_{t=1}^n\{\log(|y_t|/|y_{t-1}|)\}^2-\widehat{\gamma}_n^2$.
Under $H_0$, it is not hard to prove that $T_n\Longrightarrow N(0,
1)$ as $n\to\infty$. Thus, $H_0$ is rejected at the significance
level $\beta\in(0, 1)$ when $|T_n|>\Phi^{-1}(\beta/2)$, where
$\Phi(\cdot)$ is the cumulative distribution function of $N(0,1)$.

\section{Quasi-maximum likelihood estimation}

Let $\theta=(\phi,\alpha)^\T$ be the unknown parameter of model
(\ref{model}) with true value $\theta_{0}=(\phi_{0},\alpha_{0})^\T$.
Denote the parameter space by $\Theta=\mathcal{R}\times(0,\infty)$.
Assume that $\{y_{1},\cdots,y_{n}\}$ are generated by model
(\ref{model}). When $\eta_{t}\sim N(0,1)$, the log-likelihood
function (ignoring a constant) can be written as
\begin{eqnarray*}
L_{n}(\theta)=\sum_{t=1}^{n}l_{t}(\theta) \quad\mbox{ with }\quad
\ell_t(\theta)=-\frac{1}{2}\left\{\log(\alpha
y_{t-1}^2)+\frac{(y_t-\phi y_{t-1})^2}{\alpha y_{t-1}^2}\right\}.
\end{eqnarray*}
Then, the quasi-maximum likelihood estimator (QMLE) of $\theta_{0}$
is defined as
$$\widehat{\theta}_n:=(\widehat{\phi}_n,\widehat{\alpha}_n)^\T=\arg\max_{\theta\in\Theta}L_{n}(\theta).$$
%Let
%$\widehat{\theta}_n:=(\widehat{\phi}_n,\widehat{\alpha}_n)'=\arg\max_{\theta\in\Theta}L_{n}(\theta)$.
%Since we do not assume that $\eta_{t}\sim N(0,1)$,
%$\widehat{\theta}_n$ is called the quasi-maximum likelihood
%estimator (QMLE) of $\theta_{0}$.
By setting $\partial L_{n}(\theta)/\partial\theta=0$, it is not hard
to see that $\widehat{\theta}_n$ has a unique explicit expression
with
\begin{eqnarray*}
\widehat{\phi}_n=\frac{1}{n}\sum_{t=1}^n\frac{y_t}{y_{t-1}}\quad\mbox{
and }\quad
\widehat{\alpha}_n=\frac{1}{n}\sum_{t=1}^n\frac{(y_t-\widehat{\phi}_ny_{t-1})^2}{y_{t-1}^2}.
\end{eqnarray*}
With the help of this explicit expression, it
is convenient to obtain the asymptotic properties of
$\widehat{\theta}_n$ in the following theorem.

% Clearly, this is not the case for the QMLE of model (\ref{dar}) in
% \cite{ling04} and \cite{cll}.

{\bf Theorem 4.1.} {\it
Suppose that $E\eta_t =0$ and $E\eta_t^2=1$. Then,
\begin{description}
\item
$\mathrm{(i)}$ $\widehat{\phi}_n$ is unbiased, but
$\widehat{\alpha}_n$ is asymptotically unbiased;
\item $\mathrm{(ii)}$
$\widehat{\theta}_n\rightarrow\theta_0$ a.s. as $\rightarrow\infty$;
\item $\mathrm{(iii)}$ Furthermore, if $E\eta_{t}^{3}=0$ and $\kappa_4=E\eta_t^4<\infty$, then
$$
\sqrt{n}(\widehat{\theta}_n-\theta_0)\Longrightarrow N(0, \Sigma)
$$
as $n\rightarrow\infty$, where $\Sigma=\mathrm{diag}(\alpha_0,
(\kappa_4-1)\alpha_0^2)$.
\end{description}
}%\end{theorem}

%\begin{remark}
{\it Remark 4.1.} In Theorem 4.1 (i) and (ii),  the symmetric condition of $\eta_t$ is not required but we need one more condition that
$E\eta_{t}^{3}=0$ for part (iii) of this theorem. This is to guarantee the existence of covariance $\Sigma$. From the proof of this
theorem,  we have
$$\lim_{n\to\infty}\mathrm{cov}(\sqrt{n}(\widehat{\phi}_n-\phi_0),\sqrt{n}(\widehat{\alpha}_n-\alpha_0))=
\alpha_{0}^{3/2}E\eta_{t}^{3}\times
\left[\lim_{n\to\infty}\frac{1}{n}\sum_{t=1}^{n}E[\mbox{sign}(y_{t-1})]\right];$$
since $y_{t}$ is non-stationary, the limit in the proceeding
equation may not exist unless $E\eta_{t}^{3}=0$.%end{remark}

%\begin{remark}
{\it Remark 4.2.}
Although $\widehat{\alpha}_n$ is biased, its re-scaled version
$\widehat{\alpha}^*_n=\frac{n}{n-1}\widehat{\alpha}_n$ is unbiased.
%\end{remark}

%\begin{remark}
{\it Remark 4.3.}
Compared to the asymptotic property of the QMLE in Ling (2004) and Chen, Li and Ling (2014), Theorem 4.1 has two interesting features.
First, the compactness
of the parameter space $\Theta$ is not necessary. Second, the asymptotic covariances of $\widehat{\theta}_n$ are the
same in both the stable and unstable cases.
%\end{remark}

%\begin{remark}
{\it Remark 4.4.}
As a possible application of Theorem 4.1, one may
consider a natural plug-in estimator of $\gamma_0$ as $
\widehat{\gamma}_n^*=\frac{1}{n}\sum_{t=1}^n\log\big|\widehat{\phi}_n+\widehat{\eta}_t\sqrt{\widehat{\alpha}_n}\big|
$, where
$\widehat{\eta}_t=(y_t-\widehat{\phi}_ny_{t-1})/\sqrt{\widehat{\alpha}_ny_{t-1}^2}$
% \begin{eqnarray}\label{residual}
% \widehat{\eta}_t=\frac{y_t-\widehat{\phi}_ny_{t-1}}{\sqrt{\widehat{\alpha}_ny_{t-1}^2}}
% \end{eqnarray}
is the residual of model (\ref{model}).
However, unlike Francq and Zako\"{i}an (2012), the asymptotic normality of $\widehat{\gamma}_n^*$ requires that
$E\{(\phi_{0}+\eta_{t}\sqrt{\alpha_{0}})^{-1}\}<\infty$, which does not exist generally, even for $\eta_{t}$ being standard normal. See also Li, Guo and Li (2015) for more relevant discussions.
%It would be applicable but not simpler than $\widehat{\gamma}_{n}$ in Section 3 and, moreover, the associate theories are much complicated. See also \cite{lgl} for more relevant discussions.
%\end{remark}

%Theorem \ref{theorem2} always holds regardless of the stability of the true model (\ref{model}). And
%the asymptotic covariances of $\widehat{\theta}_n$ are always the same in both stable case and unstable ones.
%This point differ with those in \cite{ling04} and \cite{cll}. Finally, the compactness of the parameter space
%$\Theta$ is not necessary.

To construct confidence intervals for $\theta_0$, we need to
estimate $\kappa_4$. Define
$\widehat{\kappa}_n=\frac{1}{n}\sum_{t=1}^n\widehat{\eta}_t^{\,4}$.
By a simple calculation, we can show that
$\widehat{\kappa}_n\rightarrow\kappa_4$ a.s. as $n\to\infty$. Moreover, by letting
$\widehat{\Sigma}_{n}=\mathrm{diag}(\widehat{\alpha}_n,
(\widehat{\kappa}_{n}-1)\widehat{\alpha}_{n}^2)$, we can construct a
Wald test statistic
\begin{eqnarray*}
W_{n}=n(\Gamma\widehat{\theta}_{n}-r)^\T(\Gamma\widehat{\Sigma}_{n}\Gamma^\T)^{-1}(\Gamma\widehat{\theta}_{n}-r)
\end{eqnarray*}
to detect the linear null hypothesis $H_{0}:\Gamma\theta_{0}=r$,
where $\Gamma\in\mathcal{R}^{s\times2}$ is a constant matrix with
rank $s$ and $r\in\mathcal{R}^{s\times1}$ is a constant vector. At
the significance level $\beta\in(0,1)$, we reject $H_{0}$ if
$W_{n}>\Psi_{s}^{-1}(1-\beta)$, where $\Psi_{d}(\cdot)$ is the
cumulative distribution function of $\chi^2_d$. Otherwise, $H_{0}$
is not rejected.

To end this section, we offer some discussions on the model
checking. In the context of non-stationary time series, Ling {\it et al.}(2013)
considered two portmanteau tests for non-stationary ARMA models,
where one is based on the residual autocorrelation functions (ACFs)
as in Ljung and Box (1978), and the other is based on the squared residual
ACFs as in McLeod and Li (1983). However, their methods are hard to be
implemented for model (\ref{model}). To see it clearly, define the
lag-$k$ residual ACF of $\{\widehat{\eta}_t\}$ as
\begin{eqnarray*}
\breve{\rho}_k^{\,*}=\frac{\sum_{t=k+1}^n(\widehat{\eta}_t-\overline{\eta_n^*})(\widehat{\eta}_{t-k}-\overline{\eta^*_n})}
{\sum_{t=1}^n(\widehat{\eta}_t-\overline{\eta^*_n})^2},\quad
k=1,2,...,
\end{eqnarray*}
where $\overline{\eta_n^*}=\frac{1}{n}\sum_{t=1}^n\widehat{\eta}_t$.
By using the fact that
\begin{eqnarray*}
\widehat{\eta}_t=\frac{\phi_0-\widehat{\phi}_n}{\sqrt{\widehat{\alpha}_n}}\mathrm{sign}(y_{t-1})
+\eta_t\sqrt{\frac{\alpha_0}{\widehat{\alpha}_n}}\quad
\mbox{and}\quad
\lim_{n\to\infty}\frac{1}{n}\sum_{t=1}^n\widehat{\eta}_t^2\rightarrow1
\,\, \mbox{a.s.},
\end{eqnarray*}
a simple calculation entails that
\begin{eqnarray*}
\sqrt{n}\breve{\rho}_k^{\,*}=
\frac{1}{\sqrt{n}}\sum_{t=k+1}^n\eta_t\eta_{t-k}-\frac{\sqrt{n}(\widehat{\phi}_n-\phi_0)}{\sqrt{\alpha_0}}
\frac{1}{n}\sum_{t=k+1}^n\eta_{t-k}\mathrm{sign}(y_{t-1})+o_p(1).
\end{eqnarray*}
Since $\{\eta_{t-k}\mathrm{sign}(y_{t-1})\}$ is neither a stationary
nor martingale difference sequence, it is hard to determine the
limit of
$\frac{1}{n}\sum_{t=k+1}^n\eta_{t-k}\mathrm{sign}(y_{t-1})$, and
hence that of $\sqrt{n}\breve{\rho}_k^{\,*}$. Thus, the classical
portmanteau tests are not feasible for model (\ref{model}), and how
to check the adequacy of model (\ref{model}) is still a challenging
open question.

\section{Simulation studies}

In this section, we carry out simulation studies to assess the
performance of the estimator of $\gamma_{0}$, the test of stability,
and the QMLE of $\theta_{0}$ in finite samples. We generate 1000 replications of
sample size $n=100$ and 200 from the following DARWIN(1) model:
$$y_{t}=0.5y_{t-1}+\eta_{t}\sqrt{\alpha_{0}y_{t-1}^{2}},$$
where $\eta_t$ is taken as $N(0, 1)$, the standardized Student's
$t_5$ ($\mathrm{st}_5$) with density
$f(x)=\frac{8}{3\pi\sqrt{3}}(1+x^2/3)^{-3}$, and the Laplace
distribution with density
$f(x)=\frac{1}{\sqrt{2}}\exp(-\sqrt{2}|x|)$, respectively. Here, we
set $\phi_0=0.5$ and let the value of $\alpha_0$ vary corresponding
to the cases of $\gamma_{0}>0$, $\gamma_{0}=0$, and $\gamma_{0}<0$,
respectively; see Table \ref{values}.
\begin{table}[htp]
\caption{\label{values} The values of $\gamma_0$ and
$\sigma_{\gamma_{0}}^2$ when $\phi_0=0.5$ is fixed and $\alpha_{0}$
varies.}
\centering\addtolength{\tabcolsep}{-0.4pt}
\fbox{ \centering
\begin{tabular}{lrrrlrrrlrr}
 & \multicolumn{1}{c}{$N(0, 1)$} & \multicolumn{1}{c}{}&
 &  & \multicolumn{1}{c}{$\mathrm{st}_5$} & \multicolumn{1}{c}{}&
 &  & \multicolumn{1}{c}{Laplace} & \multicolumn{1}{c}{}   \\
\cline{1-3}\cline{5-7}\cline{9-11} $\alpha_0$ &
\multicolumn{1}{c}{$\gamma_0$} &
\multicolumn{1}{c}{$\sigma_{\gamma_{0}}^2$}&
 & $\alpha_0$ & \multicolumn{1}{c}{$\gamma_0$} & \multicolumn{1}{c}{$\sigma_{\gamma_{0}}^2$}&
 & $\alpha_0$ & \multicolumn{1}{c}{$\gamma_0$} & \multicolumn{1}{c}{$\sigma_{\gamma_{0}}^2$}   \\
\hline
3.1    &-0.0297  &1.2326  & &4.1 &-0.0289 &1.3355 & &5.0&-0.0143&1.4357\\
3.3058 &0.0000  &1.2328  & &4.3697 &0.0000 &1.3368 & &5.1726 &0.0000 &1.4396 \\
3.5    &0.0265  &1.2326  & &4.5 &0.0133 &1.3374 & &5.4&0.0182&1.4443\\
\end{tabular}
\hfill }
\end{table}
\begin{table}[htp]
\caption{\label{qmle} Summary for the QMLE $\widehat{\theta}_n$ and
the proposed estimator $\widehat{\gamma}_n$.}
\small\centering\addtolength{\tabcolsep}{-0.1pt}
\fbox{
\begin{tabular}{llrcrrrrrrr}
$\eta$ &$\alpha_0$ & \multicolumn{1}{c}{$\gamma_0$} &
&\multicolumn{1}{c}{} & \multicolumn{1}{c}{$n=100$} &
\multicolumn{1}{c}{}&  &\multicolumn{1}{c}{} &
\multicolumn{1}{c}{$n=200$} & \multicolumn{1}{c}{}\\
\cline{5-7}\cline{9-11}
  &  & \multicolumn{1}{c}{} & &\multicolumn{1}{c}{$\widehat{\phi}_n$} &
\multicolumn{1}{c}{$\widehat{\alpha}_n$} &
\multicolumn{1}{c}{$\widehat{\gamma}_n$}&
&\multicolumn{1}{c}{$\widehat{\phi}_n$} &
\multicolumn{1}{c}{$\widehat{\alpha}_n$} & \multicolumn{1}{c}{$\widehat{\gamma}_n$}\\
\cline{4-11}
              &3.1 &-0.0297  &EM  &0.4995  &3.0672  &-0.0308 & &0.4966&3.0845&-0.0304\\
              &    &         &ESD &0.1721  &0.4442  &0.1138  & &0.1226&0.3182& 0.0776 \\
              &    &         &ASD &0.1761  &0.4384  &0.1110  & &0.1245&0.3100& 0.0785\\
              &3.3058 &0     &EM  &0.4978  &3.2943  &-0.0029 & &0.5032&3.2912& 0.0016\\
$N(0, 1)$     &    &         &ESD &0.1767  &0.4795  &0.1155  & &0.1253&0.3291& 0.0822 \\
              &    &         &ASD &0.1818  &0.4675  &0.1110  & &0.1286&0.3306& 0.0785\\
              &3.5 &0.0265   &EM  &0.5086  &3.4548  &0.0200  & &0.4997&3.4721& 0.0247\\
              &  &           &ESD &0.1866  &0.4877  &0.1143  & &0.1292&0.3464& 0.0772 \\
              &  &           &ASD &0.1871  &0.4950  &0.1110  & &0.1323&0.3500& 0.0785\\
\hline
              &4.1 &-0.0289  &EM  &0.5013  &4.1054  &-0.0236 & &0.4985&4.1129& -0.0293\\
              &  &           &ESD &0.2030  &1.0845  &0.1181  & &0.1442&0.9438& 0.0793 \\
              &  &           &ASD &0.2025  &1.1597  &0.1156  & &0.1432&0.8200& 0.0817\\
              &4.3697 &0     &EM  &0.5162  &4.3623  &0.0037  & &0.4960&4.4137& 0.0013\\
$\mathrm{st}_5$&  &          &ESD &0.2127  &1.1575  &0.1160  & &0.1497&0.9874& 0.0796 \\
              &  &           &ASD &0.2090  &1.2359  &0.1156  & &0.1478&0.8739& 0.0818\\
              &4.5 &0.0133   &EM  &0.5066 & 4.4341 & 0.0172 & &0.4906  &4.4976  &0.0138\\
              &  &           &ESD &0.2115 & 1.1814 & 0.1235 & &0.1463  &0.8349  &0.0810\\
              &  &           &ASD &0.2121 & 1.2728 & 0.1156 & &0.1500  &0.9000  &0.0818\\
\hline
              &5.0 &-0.0143  &EM   &0.5050  &5.0073  &-0.0107 & &0.4921& 4.9650& -0.0159\\
              &    &         &ESD  &0.2229  &1.0879  &0.1184  & &0.1616& 0.7920& 0.0868 \\
              &    &         &ASD  &0.2236  &1.1180  &0.1198  & &0.1581& 0.7906& 0.0847\\
              &5.1726 &0     &EM   &0.4963  &5.1454  &0.0055  & &0.5005& 5.1590& -0.0034\\
Laplace       &  &           &ESD  &0.2241  &1.1739  &0.1195  & &0.1598& 0.8347& 0.0838 \\
              &  &           &ASD  &0.2274  &1.1566  &0.1200  & &0.1608& 0.8179& 0.0848\\
              &5.4  &0.0182  &EM   &0.4962  &5.3933  &0.0148  & &0.5023& 5.4004& 0.0195\\
              &  &           &ESD  &0.2305  &1.2174  &0.1208  & &0.1677& 0.8702& 0.0844 \\
              &  &           &ASD  &0.2324  &1.2075  &0.1202  & &0.1643& 0.8538& 0.0850\\
\end{tabular}
\hfill}
\end{table}

\begin{figure}[htp]
\begin{center}
% The arguments in the next line are {height}{optional width}{used only by OUP for typesetting}[filename, in directory art]
%\figurebox{20pc}{35pc}{}[gammalim.eps]
\includegraphics[height= 100mm,width= 140mm]{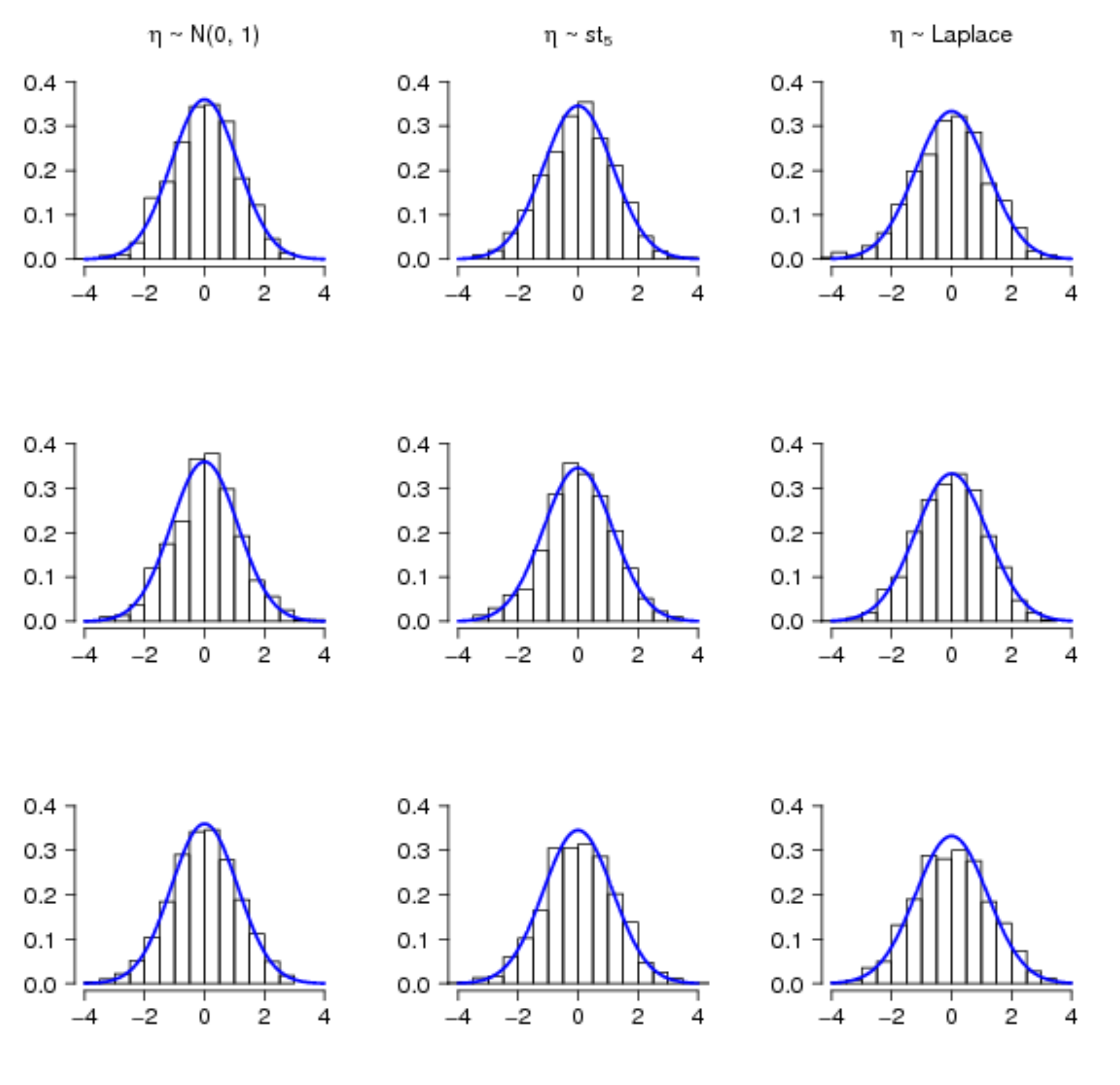}
% note that files may not be rotated
\caption{The histograms of $\sqrt{n}(\widehat{\gamma}_n-\gamma_0)$ when $\eta$ is $N(0, 1)$,
  $\mathrm{st}_5$,
   and Laplace distribution, respectively. The curves are the densities of $N(0, \sigma_{\gamma_{0}}^2)$.
   Here the true parameter $\phi_0=0.5$ and the values of $\alpha_0$ and $\sigma_{\gamma_{0}}^2$ are given in Table \ref{values},
   respectively. The top, middel, and bottom panels correspond to $\gamma_0<0$, $\gamma_0=0$, and
   $\gamma_0>0$, respectively. The sample size is 200.}
\label{gammalimdist}
\end{center}
\end{figure}

Table \ref{qmle} shows the empirical mean (EM), empirical standard
deviation (ESD), and asymptotic standard deviation (ASD) of the QMLE
$\widehat{\theta}_{n}$ and the Lyapunov exponent estimator
$\widehat{\gamma}_{n}$. The ASD of $\widehat{\gamma}_{n}$ and
$\widehat{\theta}_{n}$ is calculated by the asymptotic covariance
matrix in Theorems 3.1 and 4.1, respectively,
where the theoretical values of $\sigma_{\gamma_{0}}^{2}$ are given
in Table \ref{values}, and the theoretical values of $\kappa_4$ are
3, 9 and 6 for $N(0, 1)$, $\mathrm{st}_5$, and Laplace
distributions, respectively. From this table, we can see that the
larger the sample size, the closer the EMs and their corresponding
true values, and also the closer the ESDs and ASDs. Particularly,
$\widehat{\gamma}_n$ performs well even though $\gamma_0$ is very
small. To assess the overall performance of $\widehat{\gamma}_n$,
Fig.\ref{gammalimdist} plots the empirical density of
$\sqrt{n}(\widehat{\gamma}_n-\gamma_0)$ for different $\alpha_0$ in
Table \ref{values}. From Fig.\ref{gammalimdist}, we find that
$\widehat{\gamma}_n$ has a good performance in all cases.

\begin{figure}[htp]
\begin{center}
% The arguments in the next line are {height}{optional width}{used only by OUP for typesetting}[filename, in directory art]
%\figurebox{15pc}{35pc}{}[Tpower.eps]
\includegraphics[height= 80mm,width= 140mm]{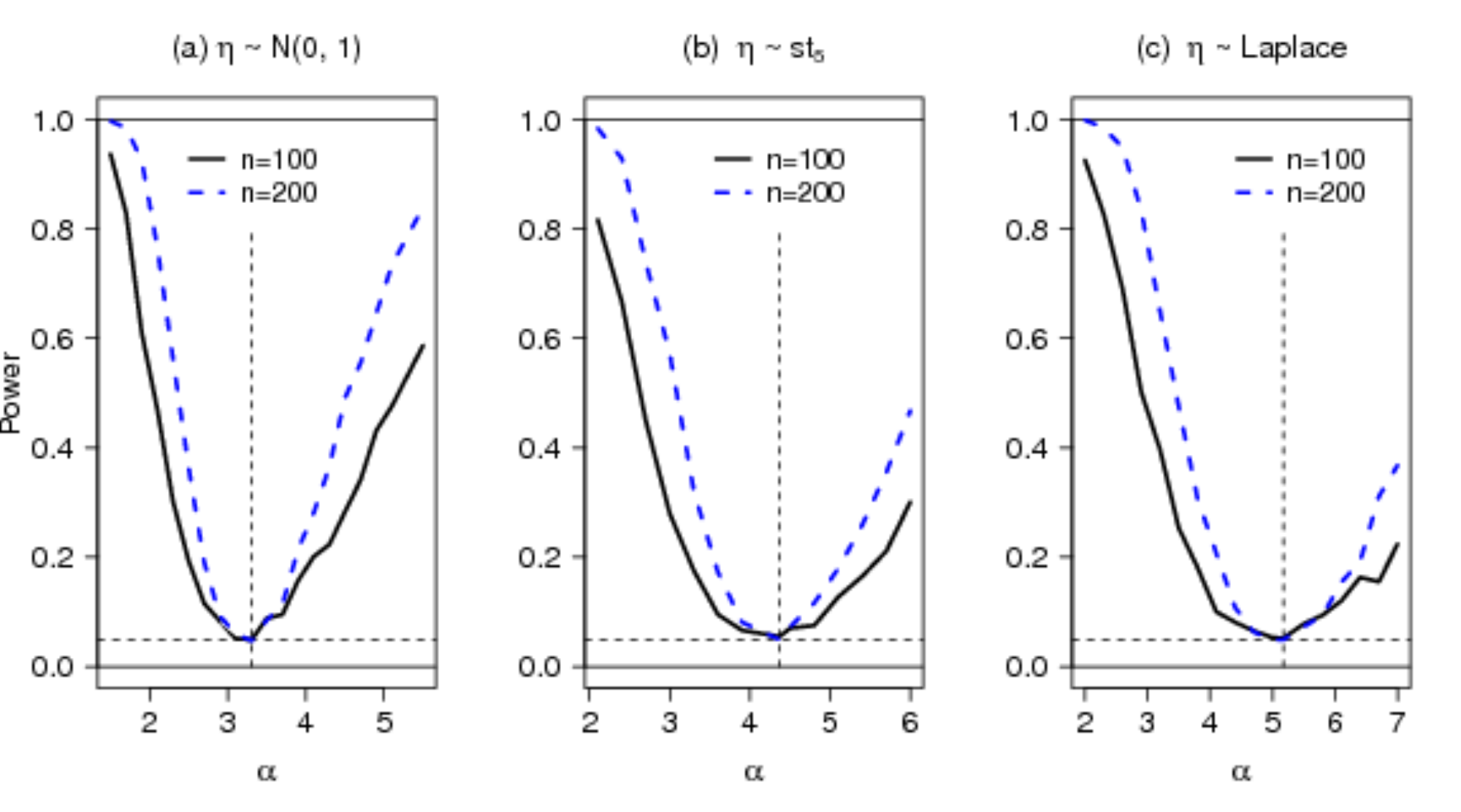}
% note that files may not be rotated
\caption{The power and size of $T_n$ at the significance level 5\% when $\eta$ is $N(0, 1)$, $\mathrm{st}_5$ and
   Laplace distribution, respectively.}
\label{testpower}
\end{center}
\end{figure}

Next, we examine the performance of test statistic $T_n$ in
(\ref{test}) for testing the hypothesis $H_0: \gamma_0=0$ against
$H_1: \gamma_0\neq 0$. Fig. \ref{testpower} shows the power and size
of $T_n$ for $n=100$ and 200 when $\alpha_{0}$ varies, and the size
of $T_{n}$ corresponds to the case that $\alpha_{0}=3.3058$ (for
$N(0, 1)$), 4.3697 (for $\mathrm{st}_5$) and 5.1726 (for Laplace);
see Table \ref{values}. From this figure, we can see that $T_n$ has
a very precise size, and overall, the power to detect the instability
is significant, even when sample size is small.

\section{An empirical example}
\begin{figure}[htp]
% The arguments in the next line are {height}{optional width}{used only by OUP for typesetting}[filename, in directory art]
%\figurebox{15pc}{35pc}{}[Exchange.eps]
\includegraphics[height= 80mm,width= 140mm]{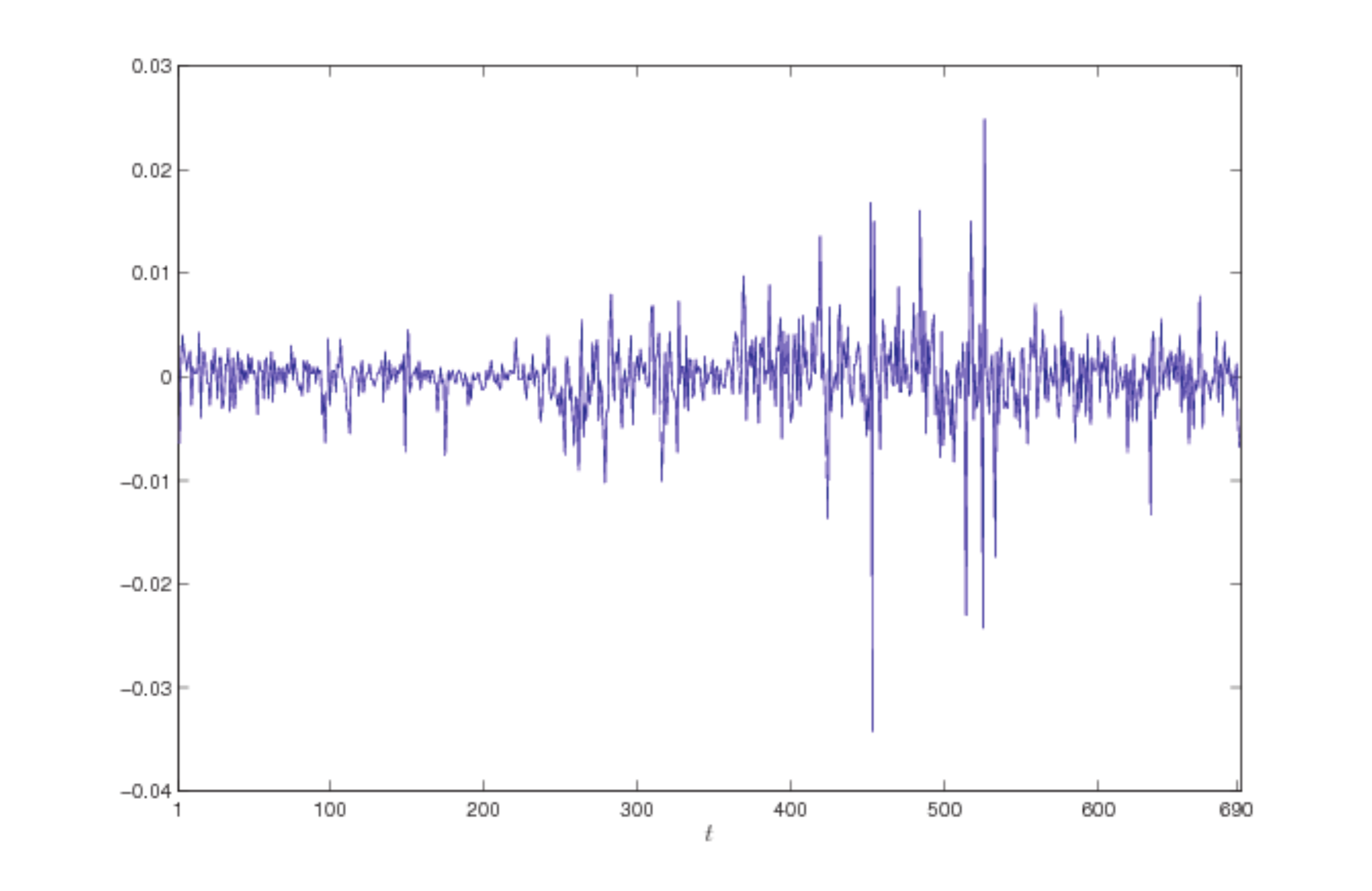}
% note that files may not be rotated
\caption{The log-return of daily exchange rates of New Taiwan Dollars
(TWD) to United States Dollars (USD) from January 1, 2007 to
December 31, 2009.}
\label{exchange_rate}
\end{figure}

This section applies the DARWIN(1) model  to study the daily
exchange rates of New Taiwan Dollars (TWD) to United States Dollars
(USD) from January 1, 2007 to December 31, 2009, which has in total
692 observations. The log-returns of this exchange rate series,
denoted by $\{y_{t}\}_{t=1}^{691}$, are plotted in
Fig.\ref{exchange_rate}.

\begin{figure}[htp]
\begin{center}
% The arguments in the next line are {height}{optional width}{used only by OUP for typesetting}[filename, in directory art]
%\figurebox{20pc}{35pc}{}[res_DARWIN.eps]
\includegraphics[height= 80mm,width= 140mm]{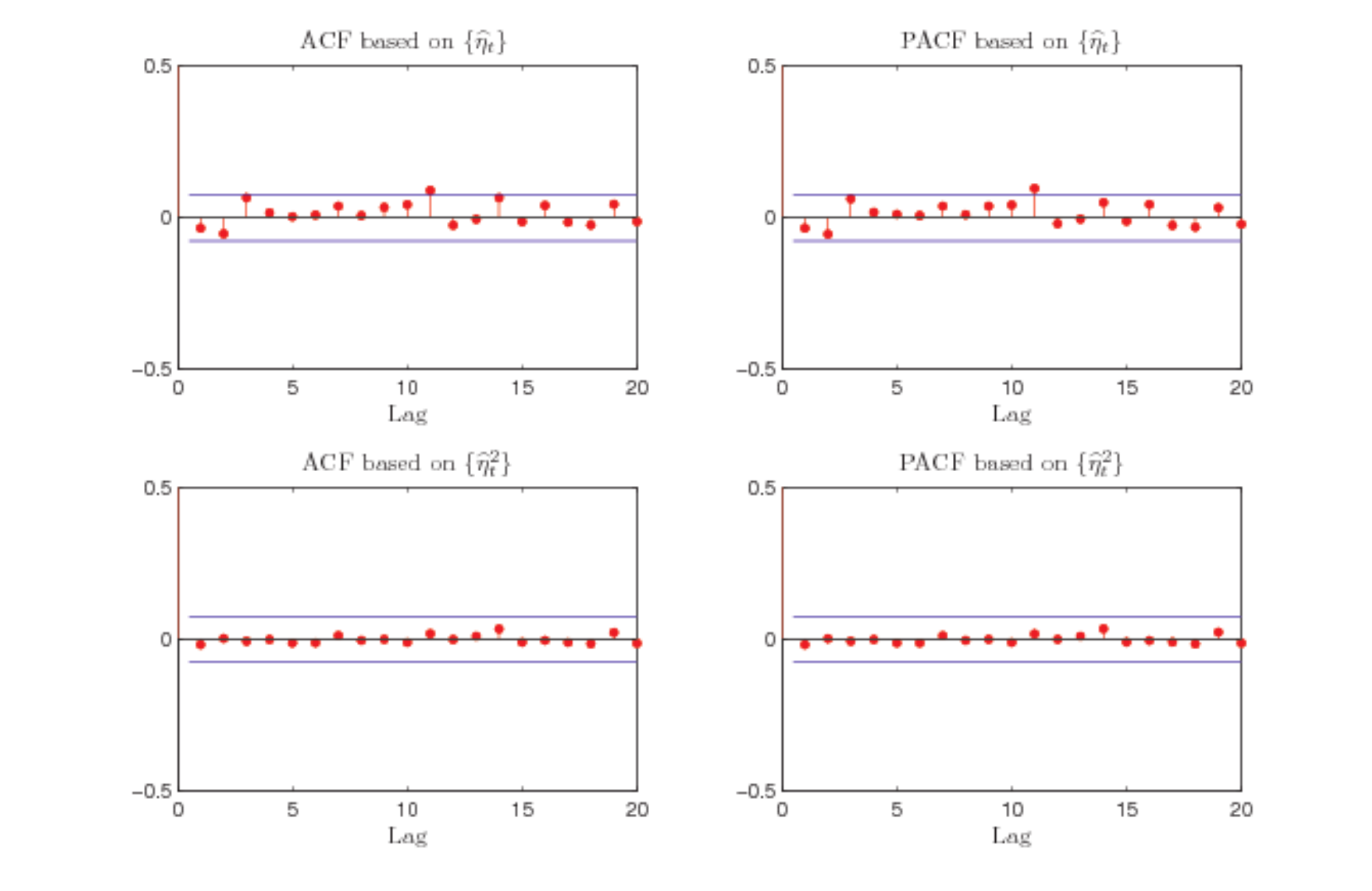}
% note that files may not be rotated
\caption{The top (or bottom) panel is the ACF and PACF of $\{\widehat{\eta}_{t}\}$ (or $\{\widehat{\eta}_{t}^{2}\})$}
\label{res_DARWIN}
\end{center}
\end{figure}

First, we use model (\ref{model}) with the QMLE estimation to fit
$\{y_t\}$ by
\begin{equation}\label{fitting}
y_t=0.3666_{(0.1661)}y_{t-1}+\eta_t\sqrt{19.0586_{(5.3083)}\,y_{t-1}^2},
\end{equation}
where the values in parentheses are estimated standard errors. Based on the residuals
$\{\widehat{\eta}_{t}\}$, Fig.\ref{res_DARWIN} plots the ACF and
PACF of $\{\widehat{\eta}_{t}\}$ and $\{\widehat{\eta}_{t}^{2}\}$.
From this figure, it seems that model (\ref{fitting}) is adequate.
Next, we use the Mira test and the Cabilio-Masaro test in \textsf{R}
package \textsf{lawstat} to test the symmetry of $\eta_{t}$, and
find that their p-values are 0.2985 and 0.2280, respectively.
Therefore, we accept the hypothesis that $\eta_{t}$ is symmetric at
the significance level 5\%. Moreover, we use the Wald test statistic
$W_{n}$ to detect the hypothesis $H_{0}: \phi_{0}=0$. The p-value of
$W_{n}$ is  0.0273, and it turns out that we can reject $H_{0}$ at
the significance level 5\%. On the other hand, we find that the
estimated Lyapunov exponent $\widehat{\gamma}_n=0.0001$ with
$\widehat{\sigma}_{\gamma}^2=1.4651$, and this implies that the
value of test statistic $T_n=0.0019$. Clearly, the null hypothesis
of $\gamma_0=0$ is not rejected at the significance level 5\%. Thus,
there is no statistical evidence against the hypothesis that the
return process is stable, i.e., log volatility of $y_{t}$ is a
random walk.

\begin{figure}[htp]
\begin{center}
% The arguments in the next line are {height}{optional width}{used only by OUP for typesetting}[filename, in directory art]
%\figurebox{20pc}{35pc}{}[res_DAR.eps]
\includegraphics[height= 80mm,width= 140mm]{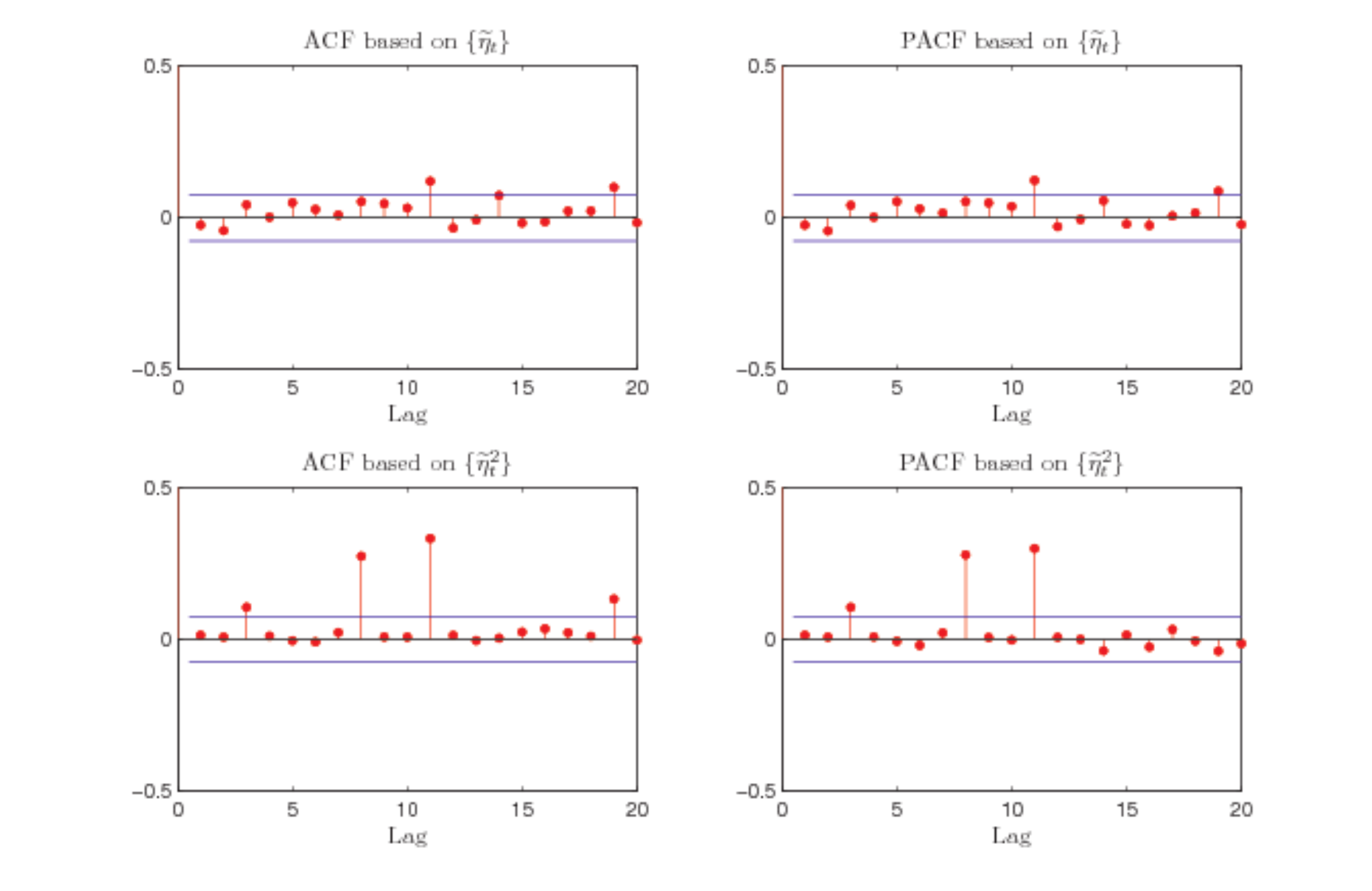}
% note that files may not be rotated
\caption{The top (or bottom) panel is the ACF and PACF of $\{\widetilde{\eta}_{t}\}$ (or $\{\widetilde{\eta}_{t}^{2}\})$}
\label{res_DAR}
\end{center}
\end{figure}

\begin{figure}[htp]
\begin{center}
% The arguments in the next line are {height}{optional width}{used only by OUP for typesetting}[filename, in directory art]
%\figurebox{20pc}{35pc}{}[log_vol.eps]
\includegraphics[height= 90mm,width= 140mm]{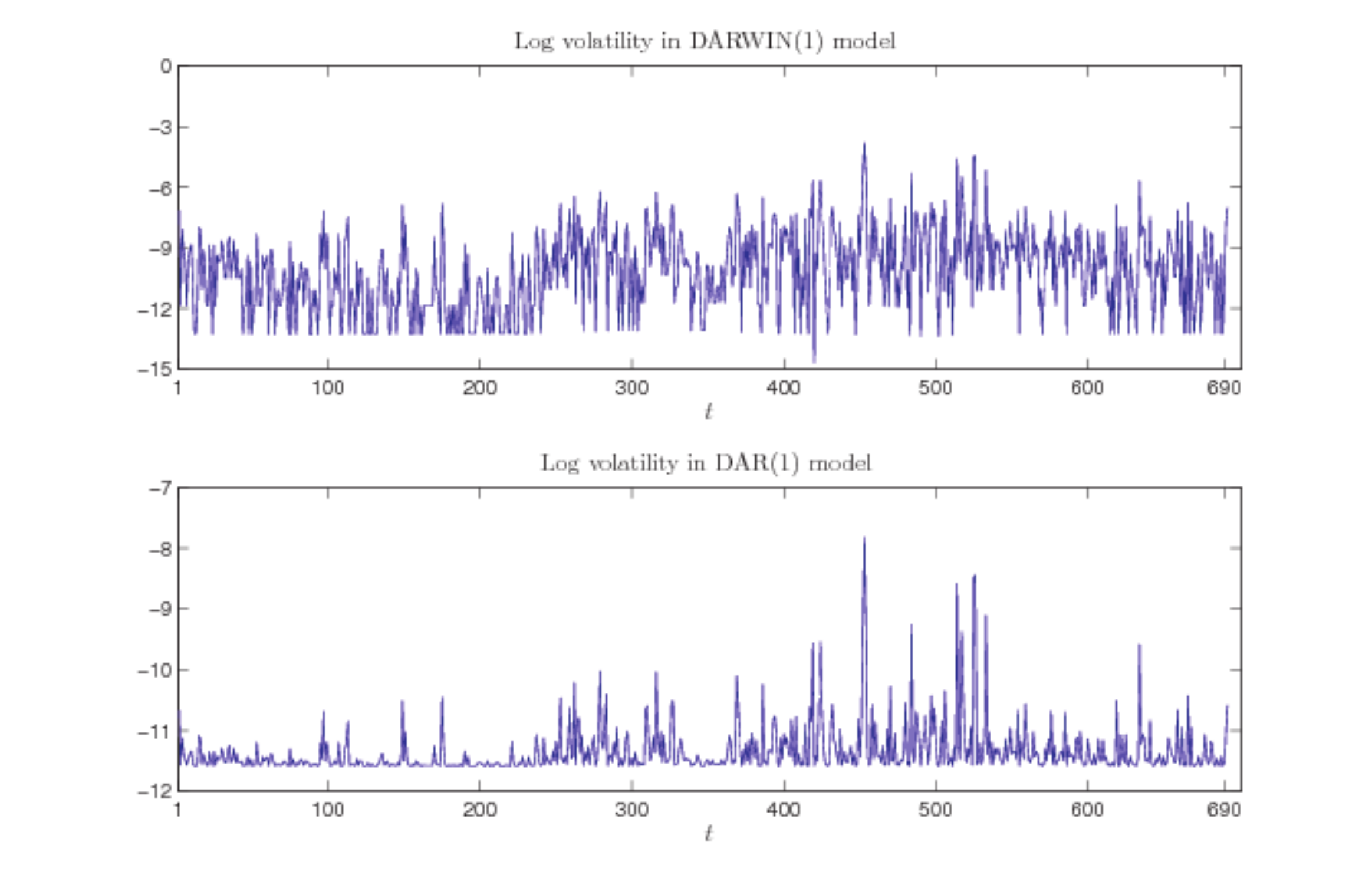}
% note that files may not be rotated
\caption{The top (or bottom) is the log volatility in model (\ref{fitting}) (or (\ref{dar_fitting})).}
\label{LogVolatility}
\end{center}
\end{figure}

It is also interesting to apply model (1) to fit $\{y_{t}\}$. The fitted model with the QMLE
estimation is
\begin{equation}\label{dar_fitting}
y_t=0.0317_{(0.0321)}
y_{t-1}+\eta_t\sqrt{9.3\times10^{-6}_{(0.1030)}+0.3349_{(0.0214)}y_{t-1}^2},
\end{equation}
where the values in parentheses are estimated standard errors. Clearly, the
estimate of $\omega_{0}$ is very close to zero. Let
$\{\widetilde{\eta}_{t}\}$ be the residuals of model
(\ref{dar_fitting}). Fig.\ref{res_DAR} plots the ACF and PACF of
$\{\widetilde{\eta}_{t}\}$ and $\{\widetilde{\eta}_{t}^{2}\}$, suggesting that model (\ref{dar_fitting}) is not
adequate. Based on these facts, it implies that model
(\ref{fitting}) is preferred to fit $\{y_{t}\}$. To gain more insight, we plot the estimated log
volatilities of models (\ref{fitting}) and (\ref{dar_fitting}) in Fig.\ref{LogVolatility}. The
log volatility of model (\ref{dar_fitting}) is bounded from below by
the logarithm of the intercept term, while there is no such bound in
model (\ref{fitting}). Among years 2007-2009, the financial crisis
happened so that the log volatilities probably tend to have no lower
bound, and hence this might lead to the preference of model
(\ref{fitting}) in fitting $\{y_{t}\}$.

%\begin{figure}[!htbp]
%  \centering
%  \includegraphics[scale=0.6]{log_vol.eps}
%  \caption{The top (or bottom) is the log volatility in model (\ref{fitting}) (or (\ref{dar_fitting})).}
%\label{LogVolatility}
%\end{figure}

\section{Concluding remarks}

This paper proposes a DARWIN(1) model. This new model is
non-stationary and heteroskedastic, but unlike non-stationary DAR model, it includes a
close-to-unit root type behavior of volatility in the stability
case. In order to study the stability of the DARWIN(1) model, an
easy-to-implement estimator of the Lyapunov exponent and a test for
the stability have been constructed. Moreover, this paper studies the QMLE of the
DARWIN(1) model, and finds that the asymptotic theory of the
QMLE is invariant regardless of the stability of the model.
An empirical example on TWD/USD exchange
rates illustrates the importance of the DARWIN(1) model.
%is preferable to fit the
%economic or financial data.

As one natural extension work, one can apply the principle of setting
the intercept to zero to high-order DAR models or conditional
heteroscedastic models, which would produce a branch of models that
totally differ from the classical counterparts in the literature.
Finally, it is worthy noting that
how to test the DARWIN model
in the null against the counterpart DAR model in the alternative is an interesting open question.
The traditional tests, e.g., the likelihood ratio test, the Lagrange
multiplier test and the Wald test, or their one-sided counterparts,
may not work due to the nonstationarity of DARWIN models. How to develop powerful tests for
such hypothesis would be a promising but challenging direction for future study.

%\section*{Acknowledgement}
%Li's work is supported in part by the Start-up Fund of Tsinghua
%University (No.553310013) and the NSFC (No.11401337).
%Guo's work is supported in part by Key
%Laboratory of RCSDS, Chinese Academy of Sciences and an EPSRC research grant in United Kingdom.
%Zhu's work is
%supported in part by the NSFC (No.11201459 and 11371354), the
%President Fund of the Academy of Mathematics and System Science,
%Chinese Academy of Sciences (grant 2014-cjrwlzx-zk), and Key
%Laboratory of RCSDS, Chinese Academy of Sciences.

%\section*{Supplementary material}
%\label{SM}
%Further material such as technical details, extended proofs, code, or additional  simulations, figures and examples may appear online, and should be briefly mentioned as Supplementary Material where appropriate.  Please submit any such content as a PDF file along with your paper, entitled `Supplementary material for Title-of-paper'.  After the acknowledgements, include a section `Supplementary material' in your paper, with the sentence `Supplementary material available at \Bka\ online includes $\ldots$', giving a brief indication of what is available.  Further instructions will be given when a paper is accepted.

%\appendix
%\appendixone
\section*{Appendix: Technical Proofs}

\subsection*{Proof of Theorem 3.1} 
By the definition of $\widehat{\gamma}_n$, we have
\begin{eqnarray*}
\widehat{\gamma}_n=\frac{1}{n}\sum_{t=1}^n\log\Big(\frac{|y_t|}{|y_{t-1}|}\Big)
=\frac{1}{n}\sum_{t=1}^n\log\big|\phi_0+\eta_t\mathrm{sign}(y_{t-1})\sqrt{\alpha_0}\big|.
\end{eqnarray*}
Since
$(\eta_1\mathrm{sign}(y_{0}),...,\eta_t\mathrm{sign}(y_{t-1}))\stackrel{d}{=}(\eta_1,...,\eta_t)$
by the induction over $t\geq 1$, we have
\begin{eqnarray*}
\widehat{\gamma}_n\stackrel{d}{=}\frac{1}{n}\sum_{t=1}^n\log\big|\phi_0+\eta_t\sqrt{\alpha_0}\big|.
\end{eqnarray*}
Thus, the result holds.

\subsection*{Proof of Theorem 4.1}
A simple calculation yields that
\begin{eqnarray*}
\widehat{\phi}_n=\phi_0+\frac{\sqrt{\alpha_0}}{n}\sum_{t=1}^n\eta_t\mathrm{sign}(y_{t-1})\quad \mbox{and}\quad
\widehat{\alpha}_n=\frac{\alpha_0}{n}\sum_{t=1}^n\eta_t^2-(\widehat{\phi}_n-\phi_0)^2.
\end{eqnarray*}
Note that $\{\eta_t\mathrm{sign}(y_{t-1})\}$ is a martingale
difference with respect to $\mathcal{F}_t=\sigma(\eta_i, i\leq t)$.
Thus, $\widehat{\phi}_n$ is unbiased. Also,
$E(\widehat{\alpha}_n)=\frac{n-1}{n}\alpha_0$. Thus, (i) holds.

By Theorem 2.19 in Hall and Hedye (1980), we have
$\frac{1}{n}\sum_{t=1}^n\eta_t\mathrm{sign}(y_{t-1})\rightarrow0$
a.s., which implies that $\widehat{\phi}_n\rightarrow\phi_0$ a.s.,
and then $\widehat{\alpha}_n\rightarrow\alpha_0$ a.s. by the strong
law of large numbers. Thus, (ii) holds.

For (iii), note that
\begin{eqnarray*}
\sqrt{n}(\widehat{\phi}_n-\phi_0)&=\frac{\sqrt{\alpha_0}}{\sqrt{n}}\sum_{t=1}^n\eta_t\mathrm{sign}(y_{t-1})\,,\\
\sqrt{n}(\widehat{\alpha}_n-\alpha_0)&=\frac{\alpha_0}{\sqrt{n}}\sum_{t=1}^n(\eta_t^2-1)
-\frac{1}{\sqrt{n}}\left(\sqrt{n}(\widehat{\phi}_n-\phi_0)\right)^2.
\end{eqnarray*}
By the Cram\'{e}r-Wold device and the martingale central limit theorem in Brown (1971), the proof of (iii) is trivial.

%\appendixtwo
%\section*{Appendix 2}
%\subsection*{Technical details}
%
%Often the appendices contain technical details of the main results.
%
%\begin{theorem}
%This is another theorem.
%\end{theorem}
%
%\appendixthree
%\section*{Appendix 3}
%
%Often the appendices contain technical details of the main results:
%\begin{equation}
%\label{C1}
%a + b = c.
%\end{equation}
%
%\begin{remark}
%This is a remark concerning equations~\eqref{A1} and \eqref{C1}.
%\end{remark}

%\bibliographystyle{biometrika}
%\bibliography{paper-ref}

%\bibliographystyle{Chicago}
\begin{center}
{\bf REFERENCES}
\end{center}

\begin{description}
\item
\textsc{Billingsley, P.} (1999). \textit{Convergence of Probability Measures}(2nd). Wiley.

\item
\textsc{Borkovec, M.} (2000). Extremal behavior of the autoregressive process with ARCH(1) errors.
\textit{Stochastic Process. Appl.} \textbf{85}, 189--207.
%
%\bibitem[\protect\citeauthoryear{Borkovec} {2001}]{bork01}
%\textsc{Borkovec, M.} (2001). Asymptotic behavior of the sample autocovariance and autocorrelation function of the AR(1) process
%with ARCH(1) errors. \textit{Bernoulli} \textbf{7}, 847--872.

\item
\textsc{Brown, B. M.} (1971).  Martingale central limit theorems.
\textit{Ann. Math. Statist.}\textbf{ 42}, 59--66.

%\bibitem[\protect\citeauthoryear{Anderson} {1959}]{anderson}
%\textsc{Anderson, T.~W.} (1959).  On asymptotic distributions of estimates of parameters of stochastic difference
%equations. \textit{Ann. Math. Statist.} \textbf{30}, 676--687.

%\bibitem[\protect\citeauthoryear{Berkes, et al.} {2009}]{bhl}
%\textsc{Berkes, I.}, \textsc{Horv\'{a}th, L.} \& \textsc{Ling, S.} (2009). Estimation in nonstationary random coefficient
%autoregressive models. \textit{J. Time Series Anal.} \textbf{30}, 395--416.

%\bibitem[\protect\citeauthoryear{Billingsley} {1999}]{bill}
%\textsc{Billingsley, P.} (1999) \textit{Convergence of Probability Measures}(2nd). Wiley.

%\bibitem[\protect\citeauthoryear{Bollerslev}{1986}]{b}
%\textsc{Bollerslev, T.} (1986)  Generalized autoregressive conditional heteroskedasticity. \textit{J. Eonometrics.}
%\textbf{31}, 307--327.

%\bibitem[\protect\citeauthoryear{Borkovec} {2000}]{bork}
%\textsc{Borkovec, M.} (2000). Extremal behavior of the autoregressive process with ARCH(1) errors.
%\textit{Stochastic Process. Appl.} \textbf{85}, 189--207.
%
%\bibitem[\protect\citeauthoryear{Borkovec} {2001}]{bork01}
%\textsc{Borkovec, M.} (2001). Asymptotic behavior of the sample autocovariance and autocorrelation function of the AR(1) process
%with ARCH(1) errors. \textit{Bernoulli} \textbf{7}, 847--872.

\item
\textsc{Borkovec, M.} and  \textsc{Kl\"{u}ppelberg, C.} (2001). The tail of the stationary distribution of
an autoregressive process with ARCH(1) errors. \textit{Ann. Appl. Prob.} \textbf{11}, 1220--1241.

\item
\textsc{Bollerslev, T}, \textsc{Chou, R.Y.} and \textsc{Kroner,
K.F.} (1992). ARCH modeling in finance: A review of the theory and
empirical evidence. \textit{J. Econometrics} \textbf{52}, 5--59.

%\bibitem[\protect\citeauthoryear{Bougerol \& Picard} {1992}]{bp}
%\textsc{Bougerol, P.} \&  \textsc{Picard, N.} (1992b). Strict stationarity of generalized autoregressive processes.
%\textit{Ann. Probab.} \textbf{20}, 1714--1730.

%\bibitem[\protect\citeauthoryear{Brown}{1971}]{brown}
%\textsc{Brown, B.~M.} (1971). Martingale  central limit theorems. \textit{Ann. Math. Statist.} \textbf{42}, 59--66.

\item
\textsc{Chan, N.~H.} and  \textsc{Peng, L.} (2005). Weighted least absolute deviation estimation for an AR(1) process
with ARCH(1) errors. \textit{Biometrika} \textbf{92}, 477--484.

\item
\textsc{Chen, M.}, \textsc{Li, D.} and \textsc{Ling, S.} (2014). Non-stationarity and quasi-maximum likelihood estimation
 on a double autoregressive model. \textit{J. Time Series Anal.} \textbf{35}, 189-202.

%\bibitem[\protect\citeauthoryear{Cline}{2007a}]{clinea}
%\textsc{Cline, D.~B.~H.} (2007). Evaluating the Lyapounov exponent and existence of moments for threshold AR-ARCH models.
%\textit{J. Time Series Anal.} \textbf{28}, 241--260.
%
%\bibitem[\protect\citeauthoryear{Cline}{2007b}]{clineb}
%\textsc{Cline, D.~B.~H.} (2007).  Regular variation of order 1 nonlinear AR-ARCH models.
%\textit{Stochastic Process. Appl.} \textbf{117}, 840--861.

\item
\textsc{Cline, D.~B.~H} and \textsc{Pu, H.~H.} (2004). Stability  and the Lyapounov exponent  of threshold  AR-ARCH models.
\textit{Ann. Appl. Prob.} \textbf{14}, 1920--1949.

\item
\textsc{Engle, R.F.} (1982). Autoregressive conditional heteroscedasticity with estimates of the variance of United Kingdom inflation.
\textit{Econometrica} \textbf{50}, 987--1007.

%\bibitem[\protect\citeauthoryear{Fan and Yao}{2003}]{fy}
%Fan, J. and Yao, Q. (2003). \textit{Nonlinear time series: nonparametric and parametric methods}.
%Springer-Verlag, New York.

%\bibitem[\protect\citeauthoryear{Francq, Olivier and Zako\"{i}an}{2013}]{foz}
%Francq, C., Olivier, W. and Zako\"{i}an, J.-M. (2013). GARCH models without positivity constraints: Exponential or
%log GARCH?  \textit{J. Econometrics} \textbf{177}, 34--46.
%
%\bibitem[\protect\citeauthoryear{Francq and Zako\"{i}an}{2000}]{fz00}
%Francq, C. and Zako\"{i}an, J.-M. (2000).  Estimating weak GARCH representations.
%\textit{Econometric Theory} \textbf{26}, 692--728.
%
%\bibitem[\protect\citeauthoryear{Francq and Zako\"{i}an}{2004}]{fz04}
%Francq, C. and Zako\"{i}an, J.-M. (2004). Maximum likelihood estimation of pure GARCH and ARMA-GARCH processes.
%\textit{Bernoulli} \textbf{10}, 605--637.

%\bibitem[\protect\citeauthoryear{Francq and Zako\"{i}an}{2006}]{fz06}
%Francq, C. and Zako\"{i}an, J.-M. (2006). Mixing properties of a general class of GARCH(1,1) models
%without moment assumptions on the observed process. \textit{Econometric Theory} \textbf{22}, 815--834.
%
%\bibitem[\protect\citeauthoryear{Francq and Zako\"{i}an}{2007}]{fz07}
%Francq, C. and Zako\"{i}an, J.-M. (2007). Quasi-maximum likelihood estimation in GARCH processes
%when some coefficients are equal to zero. \textit{Stochastic Process. Appl.} \textbf{117}, 1265--1284.
%
\item
\textsc{Francq, C.} and \textsc{Zako\"{i}an, J.-M.} (2010).  \textit{GARCH Models: Structure, Statistical Inference and
Financial Applications}. John Wiley.
%
%\bibitem[\protect\citeauthoryear{Francq, Lepage and Zako\"{i}an}{2011}]{flz}
%Francq, C., Lepage, G. and Zako\"{i}an, J.-M. (2011). Two-stage non Gaussian QML estimation of GARCH models
%and testing the efficiency of the Gaussian QMLE. \textit{J. Econometrics} \textbf{165}, 246--257.

\item
\textsc{Francq, C.} and \textsc{Zako\"{i}an, J.-M.} (2012).  Strict stationarity testing and estimation of explosive
and stationary generalized autoregressive conditional heteroscedasticity models. \textit{Econometrica}  \textbf{80}, 821--861.

%\bibitem[\protect\citeauthoryear{Francq \& Zako\"{i}an}{2013}]{fz13}
%\textsc{Francq, C.} \& \textsc{Zako\"{i}an, J.-M.} (2013). Inference in nonstationary asymmetric
%GARCH models. \textit{Ann. Statist.} \textbf{41}, 1970--1998.

%\bibitem[\protect\citeauthoryear{Gouri\'{e}roux}{1997}]{g}
%Gouri\'{e}roux, C. (1997). \textit{ARCH Models and Financial Applications}. Springer-Verlag, New York.

\item
\textsc{Gu\'{e}gan, D.} and \textsc{Diebolt, J.}  (1994). Probabilistic properties of the $\beta$-ARCH-model.
 \textit{Statist. Sinica}  \textbf{4}, 71--87.

%\bibitem[\protect\citeauthoryear{Hall and Yao}{2003}]{hy}
%Hall, P. and Yao, Q. (2003). Inference in ARCH and GARCH models with heavy-tailed errors.
%\textit{Econometrica} \textbf{71}, 285--317.

%\bibitem[\protect\citeauthoryear{Hitczenko and Weso{\l}owski}{2011}]{hw}
%Hitczenko, P. and Weso{\l}owski, J. (2011). Renorming divergent perpetuities.
%\textit{Bernoulli} \textbf{17}, 880--894.

\item
Guo, S.,  Ling, S. and Zhu, K. (2014).  Factor double autoregressive models with application to simultaneous causality testing.
\textit{J. Statist. Plann. Inference}  \textbf{148}, 82--94.

%\bibitem[\protect\citeauthoryear{Hjort and Jones}{1996}]{hj}
%\textsc{Hjort, N.~L.} and \textsc{Jones, M.~C.} (1996). Better Rules of Thumb for Choosing Bandwidth in Density Estimation.
%\textit{Statistical research report}, Department of Mathematics, University of Oslo, Oslo, Norway.

\item
\textsc{Hafner, C.~M.} and \textsc{Preminger, A.} (2015). An ARCH model without intercept.
\textit{Econom. Lett.} \textbf{129}, 13--17.

\item
\textsc{Hall, P.} and \textsc{Heyde, C. C.} (1980).
\textit{Martingale limit theory and its application.}  Academic Press, New York.

%\bibitem[\protect\citeauthoryear{Hjort and Pollard}{1993}]{hp}
%Hjort, N.~L. and Pollard, D. (1993). Asymptotics for minimisers of convex processes
%\textit{Unpublished manuscrip}.

%\bibitem[\protect\citeauthoryear{Jensen and Rahbek}{2004a}]{jra}
%\textsc{Jensen, S. T.} and \textsc{Rahbek, A.} (2004a).  Asymptotic normality of the QMLE estimator of ARCH
%in the nonstationary case. \textit{Econometrica} \textbf{72}, 641--646.
%
%\bibitem[\protect\citeauthoryear{Jensen and Rahbek}{2004b}]{jrb}
%\textsc{Jensen, S. T.} and \textsc{Rahbek, A.} (2004b).  Asymptotic inference for nonstationary GARCH.
%\textit{Econometric Theory} \textbf{20}, 1203--1226.

\item
\textsc{Li, D.}, \textsc{Guo, S.} and \textsc{Li, M.} (2015). Robust
inference in a double AR model. \textit{Working paper}.

\item
\textsc{Li, D.}, \textsc{Li, M.} and \textsc{Wu, W.} (2014). On dynamics of volatilities in nonstationary GARCH models.
\textit{Statist. Probab. Lett.} \textbf{94}, 86--90.

\item
\textsc{Li, D.}, \textsc{Ling, S.} and \textsc{Zako\"{i}an, J.-M.}
(2015). Asymptotic inference in multiple-threshold double
autoregressive models. Forthcoming in \textit{J. Econometrics}.

\item
\textsc{Li, D.}, \textsc{Ling, S.} and \textsc{Zhang, R.M.} (2015).
On a threshold double autoregressive model. Forthcoming in \textit{J.
Bus. Econom. Statist.}

%\bibitem[\protect\citeauthoryear{Li}{2004}]{libook}
%\textsc{Li, W. K.} (2004). \textit{Diagnostic checks in time series}, Chapman \& Hall/CRC.

%\bibitem[\protect\citeauthoryear{Li and Mak}{1994}]{limak}
%\textsc{Li, W. K.} and \textsc{Mak, T. K.} (1994). On the squared
%residual autocorrelations in non-linear time with conditional
%heteroscedasticity. \textit{J. Time Ser. Anal.} \textbf{15},
%627--636.

\item
\textsc{Ling, S.} (2004). Estimation and testing  stationarity  for double autoregressive models.
\textit{J. Roy. Statist. Soc. Ser. B} \textbf{66}, 63--78.

\item
\textsc{Ling, S.} (2007).  A double AR(p) model: structure and estimation. \textit{Statist. Sinica}  \textbf{17}, 161--175.

\item
\textsc{Ling, S.} and \textsc{Li, D.} (2008). Asymptotic inference for a nonstationary double AR(1) model.
\textit{Biometrika} \textbf{95}, 257--263.

\item
\textsc{Ling, S.},  \textsc{Zhu, K.} and \textsc{Chong, C.Y.}
(2013). Diagnostic checking for non-stationary ARMA models with an
application to financial data. \textit{North American Journal of
Economics and Finance} \textbf{26}, 624--639.

\item
\textsc{Ljung, G. M.} and \textsc{Box, G. E. P.} (1978). On a measure of lack of fit in time series models.
\textit{Biometrika} \textbf{65}, 297--303.

\item
\textsc{Lu, Z.} (1998).  On the geometric ergodicity of a non-linear autoregressive model  with an autoregressive
conditional heteroscedastic term. \textit{Statist. Sinica} \textbf{8}, 1205--1217.

%\bibitem[\protect\citeauthoryear{Nelson}{1990}]{nelson}
%\textsc{Nelson, D. B.} (1990). Stationarity and persistence in the GARCH(1,1) model.
%\textit{Econometric Theory} \textbf{6}, 318--334.
%
%\bibitem[\protect\citeauthoryear{Pollard}{1991}]{pollard}
%\textsc{Pollard, D.} (1991). Asymptotics for least absolute deviation regression estimators.
%\textit{Econometric Theory} \textbf{7}, 186--199.

\item
\textsc{McLeod, A. I.} and \textsc{Li, W. K.} (1983). Disgnostic
checking ARMA times series models using squared-residual
autocorrelations. \textit{J. Time Ser. Anal.} \textbf{4}, 269--273.

%\bibitem[\protect\citeauthoryear{Phillips and Perron}{1988}]{pp}
%\textsc{Phillips, P.C.B.} and \textsc{Perron, P.} (1988). Testing
%for a unit root in time series regression, \textit{Biometrika}
%\textbf{75}, 335--346.

%\bibitem[\protect\citeauthoryear{Tsay}{1987}]{tsay}
%\textsc{Tsay, R. S.} (1987). Conditional heteroscedastic time series models, \textit{J. Amer. Statist. Assoc.} \textbf{82}, 590--604.

\item
\textsc{Weiss, A.~A.} (1984). ARMA models with ARCH errors. \textit{J. Time Ser. Anal.}  \textbf{5}, 129--143.

\item
\textsc{Zhu, K.} and \textsc{Ling, S.} (2013).  Quasi-maximum
exponential likelihood estimators for a double AR($p$) model.
\textit{Statist. Sinica} \textbf{23}, 251--270.

\end{description}

\end{document}